\documentclass[12pt,leqno]{article}
\usepackage{amsmath,amssymb,amsthm,amsfonts,mathrsfs}
\usepackage{appendix}
\usepackage{array,enumitem,graphics,xcolor,yfonts,colonequals}
\usepackage{stmaryrd}
\usepackage[alphabetic]{amsrefs}
\usepackage[all,cmtip]{xy}
\usepackage{tikz-cd}
\usepackage[colorlinks,anchorcolor=blue,citecolor=blue,linkcolor=blue,urlcolor=blue,bookmarksopen=true]{hyperref}
\usepackage{comment}
\usepackage{mathtools}
\usepackage[margin=1in]{geometry}

\AtBeginDocument{%
	\def\MR#1{}
}

\newcommand{\bC}{\mathbb{C}}    
\newcommand{\bP}{\mathbb{P}}    
\newcommand{\bR}{\mathbb{R}}    
\newcommand{\bZ}{\mathbb{Z}}    

\newcommand{\bS}{\mathbf{S}}    

\newcommand{\cA}{\mathcal{A}}   
\newcommand{\cD}{\mathcal{D}}
\newcommand{\cO}{\mathcal{O}}   

\newcommand{\sA}{\mathscr{A}}   
\newcommand{\sP}{\mathscr{P}}   

\newcommand{\Coh}{\mathrm{Coh}} 
\newcommand{\Cone}{\mathrm{Cone}}
\newcommand{\Db}{\mathrm{D^b}}  
\newcommand{\id}{\mathrm{id}}
\newcommand{\pol}{\mathrm{pol}}   

\newcommand{\std}{\mathrm{std}} 

\newcommand{\Aut}{\operatorname{Aut}}   
\newcommand{\Hom}{\operatorname{Hom}}   
\newcommand{\Pic}{\operatorname{Pic}}	

\newcommand{\RHom}{\mathbf{R}\text{Hom}}

\newtheorem*{thm*}{Theorem}
\newtheorem*{prop*}{Proposition}
\newtheorem*{cor*}{Corollary}
\newtheorem*{notn*}{Notation}
\newtheorem{thm}{Theorem}[section]
\newtheorem{prop}[thm]{Proposition}
\newtheorem{cor}[thm]{Corollary}
\newtheorem{lemma}[thm]{Lemma}

\numberwithin{equation}{section}

\theoremstyle{definition}
\newtheorem{defn}[thm]{Definition}
\newtheorem{eg}[thm]{Example}

\newtheorem{rmk}[thm]{Remark}
\newtheorem{notn}[thm]{Notation}

\newcounter{nameOfYourChoice}

\begin{document}
\title{On entropy of $\bP$-twists}
\author{Yu-Wei Fan}
\date{}
\maketitle

\begin{abstract}
We show that the $\bP$-twist associated to any $\bP$-object of a smooth project variety is not conjugate to a standard autoequivalence. This result is obtained by computing the categorical entropy functions of $\bP$-twists. We also determine the categorical polynomial entropy of spherical twists and $\bP$-twists, under an additional assumption. 
\end{abstract}

\section{Introduction}
Let $X$ be a smooth projective variety over $\bC$, and let $\Db(X)$ denote the bounded derived category of coherent sheaves on $X$. The group of autoequivalences $\Aut\Db(X)$ has been extensively studied in the literature. A classical result by Bondal and Orlov \cite{BonOrl} states that if the (anti-)canonical bundle $\omega_X$ is ample, then the group of autoequivalences consists of only the \emph{standard} autoequivalences:
$$
\Aut_\std\Db(X)\coloneqq\left(\Pic(X)\rtimes\Aut(X)\right)\times\bZ[1],
$$
which are generated by the automorphisms of $X$, tensoring line bundles, and shifts.

In general, $\Db(X)$ may admit non-standard autoequivalences.
For instance, when $X$ is Calabi--Yau, its derived category admits \emph{spherical objects}, which induce \emph{spherical twists} in $\Aut\Db(X)$ \cite{SeidelThomas}.
Spherical twists can be regarded as the \emph{mirror} of Dehn twists along Lagrangian spheres under mirror symmetry, and thus have their origin in symplectic geometry, unlike the standard autoequivalences.
It is proved in \cite{LiLiu}*{Theorem~A.2}, using results on \emph{categorical entropy functions} of spherical twists in \cite{OuchiSpherical}*{Theorem~1.4}, that spherical twists are not conjugate to any standard autoequivalence, which answers a question posed in \cite{Huy16}*{Chapter 16}.

The categorical entropy function $h_t(\Phi)\colon\bR\rightarrow\bR$, introduced in \cite{DHKK}, is a \emph{dynamical} invariant associated with an endofunctor of a triangulated category. Roughly speaking, it measures the exponential growth rate of certain quantities under large iterations $\Phi^n$. Its precise definition is provided in Section~\ref{subsec:DynamicalInvariants}.
Note that the entropy function is conjugacy-invariant. Therefore, to prove that spherical twists are not conjugate to any standard autoequivalence, it suffices to show that their categorical entropy functions are different.
The proof consists of the following steps:
\begin{enumerate}[label=(\roman*)]
    \item\label{step:Ouchi} The categorical entropy function of the spherical twist $T_E$ of a spherical object $E$ is computed in \cite{OuchiSpherical}*{Theorem~1.4}, under the condition that $E^\perp\neq\{0\}$ in $\Db(X)$.
    \item\label{step:LiLiu} It is proved in \cite{LiLiu}*{Proposition~1.1} that $E^\perp\neq\{0\}$ for any spherical object $E\in\Db(X)$. Thus, the calculation of the categorical entropy function of $T_E$ is complete.
    \item It is not difficult to compute the categorical entropy functions of standard autoequivalences, and show that they differ from the categorical entropy functions of spherical twists obtained in Steps~\ref{step:Ouchi} and \ref{step:LiLiu}.
\end{enumerate}

In \cite{HuyThoPtwist}, the notion of \emph{$\bP$-twists} associated to $\bP$-objects is introduced. Similar to spherical twists, $\bP$-twists have their origins in symplectic geometry, where they can be regarded as the mirror of Dehn twists along Lagrangian complex projective space.
One would expect that $\bP$-twists also are not conjugate to any standard autoequivalence.
Indeed, it is not difficult to obtain the analogous result of Step~\ref{step:Ouchi} for $\bP$-twists, see Section~\ref{subsec:P-twistOuchiVersion}.
However, verifying the condition $E^\perp\neq\{0\}$ in $\Db(X)$, as outlined in Step~\ref{step:LiLiu}, is quite challenging for a general $\bP$-object. 

In the present article, we take a different approach. By utilizing recent results on \emph{shifting numbers} of autoequivalences \cites{FanFilip,FanShifting}, we are able to compute the entropy functions of spherical twists and $\bP$-twists without assuming $E^\perp\neq\{0\}$ in $\Db(X)$. Consequently, we bypass Step~\ref{step:LiLiu} and prove the following theorem.

\begin{thm}
\label{thm:MainP-twist}
Consider a complex smooth projective variety $X$ of dimension $2d$, and let $E\in\Db(X)$ be a $\bP^d$-object with the associated $\bP^d$-twist $P_E\in\Aut\Db(X)$.
Let $k$ be a nonzero integer, and $Y$ be a smooth projective variety with an exact equivalence $\Psi\colon\Db(X)\rightarrow\Db(Y)$. 

There does not exist a standard autoequivalence $\Phi\in\Aut_\std\Db(Y)$ such that $P_E^k=\Psi^{-1}\circ \Phi\circ\Psi$. In particular, $P_E^k$ is not conjugate to any standard autoequivalence of $X$.
\end{thm}

We now describe in more detail how to remove the assumption $E^\perp\neq\{0\}$.
First, for any endofunctor $\Phi$, its entropy function $h_t(\Phi)$ grows linearly as $t\rightarrow+\infty$ (resp.~$t\rightarrow-\infty$), with slope $\tau^+(\Phi)$ (resp.~$\tau^-(\Phi)$), called the \emph{upper} (resp.~\emph{lower}) \emph{shifting number} of $\Phi$ (see Proposition~\ref{prop:ShiftingNumStrip}). Roughly speaking, shifting numbers measure the asymptotic amount by which $\Phi$ translates in the triangulated category \cite{FanFilip},\cite{FanShifting}. The shifting numbers are conjugacy invariants and always satisfy $\tau^+(\Phi)\geq\tau^-(\Phi)$. The proof of Theorem~\ref{thm:MainP-twist} consists of the following steps:
\begin{enumerate}[label=(\roman*)]
    \item For standard autoequivalences, the upper and lower shifting numbers coincide. Thus, if the strict inequality $\tau^+(P_E)>\tau^-(P_E)$ holds (which implies $\tau^+(P_E^k)>\tau^-(P_E^k)$ for any $k\neq0$), then $P_E^k$ is not conjugate to any standard autoequivalences.
    \item It is not difficult to obtain $\tau^-(P_E)=-2d$ and $\tau^+(P_E)\leq0$. The main difficulty lies in proving $\tau^+(P_E)=0$.
    \item $\tau^+(P_E)=0$ can be ensured if $E^\perp\neq\{0\}$. However, without assuming $E^\perp\neq\{0\}$, we can use results from \cite{FanShifting} to show that $\tau^+(P_E)=0$ if there exists a nonzero object $F\in\Db(X)$ and a heart $\sA\subseteq\Db(X)$ such that 
    $$
    \lim_{n\rightarrow\infty}\frac{\phi_{\sA}^+(P_E^n(F))}{n}=0,
    $$
    where $\phi_\sA^+(P_E^n(F))$ denotes the maximal degree of cohomology objects of $P_E^n(F)$ with respect to the heart $\sA$ (see Notation~\ref{notn:HeartCohomologyObjectPhase}).
    \item The desired limit is satisfied by $\cA=\Coh(X)$ and $F=\cO(-m)$ for sufficiently large $m$. This proves that $\tau^+(P_E)=0>-2d=\tau^-(P_E)$.
\end{enumerate}

We also study the \emph{categorical polynomial entropy functions} of spherical twists and $\bP$-twists. The categorical polynomial entropy function, as defined in \cite{FanFuOuchi}, is a refined secondary invariant that measures the \emph{polynomial growth}, rather than exponential growth, of certain quantities under large iterations of an endofunctor. It is proved in \cite{FanFuOuchi}*{Section~6}, that for a spherical twist (or a $\bP$-twist) along a spherical object (or a $\bP$-object) $E$, the polynomial entropy function $h_t^\pol$ is given by:
\begin{itemize}
    \item $h_t^\pol=0$ for $t<0$.
    \item $h_t^\pol=0$ for $t>0$, under the assumption that $E^\perp\neq\{0\}$.
    \item $0\leq h_0^\pol\leq 1$.
\end{itemize}
Using a similar approach as in the proof of Theorem~\ref{thm:MainP-twist}, we remove the assumption $E^\perp\neq\{0\}$ for the case $t>0$. Furthermore, we establish $h_0^\pol=1$ under an additional assumption.

\begin{thm}
\label{thm:PolyCat}
Let $X$ be a complex smooth projective variety of dimension $d\geq2$. Let $E\in\Db(X)$ be a spherical object (resp.~$\bP$-object) with the associated spherical twist $T_E\in\Aut\Db(X)$ (resp.~$\bP$-twist $P_E\in\Aut\Db(X)$). Then
\begin{enumerate}[label=(\roman*)]
    \item For $t\neq0$, $h_t^\pol(T_E)=0$ (resp.~$h_t^\pol(P_E)=0$).
    \item If there exists a Bridgeland stability condition $\sigma=(Z,\sP)$ on $\cD$ such that $E\in\sP(-1,1]$, then $h_0^\pol(T_E)=1$ (resp.~$h_0^\pol(P_E)=1$).
\end{enumerate}
\end{thm}

For instance, this fully determines the polynomial entropy function of spherical twists for K3 surfaces of Picard number one, cf.~\cite{BB17}*{Corollary~6.9}.

\ \\
\noindent\textbf{Convention.}
Throughout this article, all triangulated categories (usually denoted by $\cD$) are assumed to be $\bZ$-graded, linear over a base field $\mathbf{k}$, saturated (i.e.~admit a dg-enhancement which is smooth and proper), and of finite type (i.e.~$\oplus_{k\in\bZ}\Hom_\cD(E,F[k])$ is finite-dimensional for any pair of objects $E,F$ in $\cD$). Furthermore, we assume that $\cD$ admits a Serre functor $\bS$, and contains a split generator $G\in\cD$. Endofunctors of triangulated categories are assumed to be $\mathbf{k}$-linear, triangulated, and not virtually zero (i.e.~any power is not the zero functor).

\ \\
\noindent\textbf{Acknowledgment.}
The author would like to thank Chunyi~Li for explaining Lemma~\ref{lemma:concentrateInd-1}, and Yijia~Liu for proposing this problem.

\section{Preliminaries}

\subsection{Dynamical invariants of endofunctors}
\label{subsec:DynamicalInvariants}
In this subsection, we recall the definitions and basic properties of various dynamical invariants associated with endofunctors of triangulated categories.

Let $\Phi\colon\cD\rightarrow\cD$ be an endofunctor of a triangulated category $\cD$. The \emph{categorical entropy function} of $\Phi$, introduced in \cite{DHKK}*{Definition~2.5}, is a function $h_t(\Phi)\colon\bR\rightarrow[-\infty,\infty)$ that depends on the variable $t$.
According to \cite{DHKK}*{Theorem~2.7}, it can be expressed as follows:
\begin{equation}\label{eqn:EntropyDefn}
h_t(\Phi)=\lim_{n\rightarrow\infty}\frac{1}{n}\log\left(\sum_{k\in\bZ}\dim\Hom\left(G,\Phi^nG'[k]\right)e^{-kt}\right),
\end{equation}
where $G$ and $G'$ are split generators of $\cD$.
For convenience, we denote
$$
\epsilon_t(M,N)\coloneqq\sum_{k\in\bZ}\dim\Hom(M,N[k])e^{-kt}.
$$
We collect some of the basic properties of the categorical entropy function $h_t(\Phi)$ in the following proposition.

\begin{prop}
\label{prop:EntropyBasic}
Consider an endofunctor $\Phi\colon\cD\rightarrow\cD$ of a triangulated category $\cD$.
\begin{enumerate}[label=(\roman*)]
    \item The limit in (\ref{eqn:EntropyDefn}) exists and is independent of the choice of the generators $G,G'$. See \cite{DHKK}*{Lemma~2.6 and Theorem~2.7}.
    \item The categorical entropy function $h_t(\Phi)\colon\bR\rightarrow\bR$ is real-valued and convex. See \cite{FanFilip}*{Theorem~2.1.6}.
    \item $h_0(\Phi)\geq0$. See \cite{DHKK}*{Definition~2.5 and Theorem~2.7}.
    \item\label{propEntropy:conjugacyInvar} Let $\Psi\colon\cD\rightarrow\cD'$ be an exact equivalence. Then $h_t(\Phi)=h_t(\Psi\circ\Phi\circ\Psi^{-1})$. See \cite{Kikuta}*{Lemma~2.9}.
    \item $h_t(\Phi[\ell])=h_t(\Phi)+\ell t$ for any integer $\ell$.
    \item $h_t(\Phi^k)=kh_t(\Phi)$ for $k\in\bZ_{\geq1}$.
    \item $h_t(\Phi^{-1})=h_{-t}(\Phi)$ for an autoequivalence $\Phi$. See \cite{FanFuOuchi}*{Lemma~2.11}.
\end{enumerate}
\end{prop}

$$
\sim\sim\sim
$$

There are also invariants called the \emph{upper} and \emph{lower shifting numbers} of $\Phi$, which are related to the entropy function $h_t(\Phi)$ in the following way:

\begin{prop}[\cite{ElaginLunts}*{Proposition~6.13}, \cite{FanFilip}*{Theorem~2.1.7}]\label{prop:ShiftingNumStrip}
There exist real numbers $\tau^+(\Phi)$ and $\tau^-(\Phi)$, called the \emph{upper} and \emph{lower shifting numbers} of $\Phi$, such that $h_t(\Phi)$ is bounded within the following ranges:
\begin{align*}
    t\cdot\tau^+(\Phi)\leq h_t(\Phi)\leq t\cdot\tau^+(\Phi)+h_0(\Phi) & \qquad\text{ for } t\geq0, \\
    t\cdot\tau^-(\Phi)\leq h_t(\Phi)\leq t\cdot\tau^-(\Phi)+h_0(\Phi) & \qquad\text{ for } t\leq0.
\end{align*}
Note that $\tau^+(\Phi)\geq\tau^-(\Phi)$.
\end{prop}

\begin{rmk}
By Proposition~\ref{prop:ShiftingNumStrip}, when $h_0(\Phi)=0$ (as is the case for spherical twists and $\bP$-twists, which we will discuss later), the entropy function $h_t(\Phi)$ is completely determined by the shifting numbers $\tau^\pm(\Phi)$.
\end{rmk}

\begin{eg}
\label{eg:stdautoeq}
Let $X$ be a smooth projective variety. A standard autoequivalence $\Phi\in\Aut_{\std}\Db(X)$ is of the form
$$
\Phi=f^*\circ(-\otimes L)[\ell]
$$
for some automorphism $f\in\Aut(X)$, line bundle $L$, and integer $\ell$. By \cite{DHKK}*{Lemma~2.11}, $h_t\left(f^*\circ(-\otimes L)\right)$ is a constant function of $t$. Therefore, the entropy function of $\Phi$ is
$$
h_t(\Phi)=h_0\left(f^*\circ(-\otimes L)\right)+\ell t=h_0(\Phi)+\ell t.
$$
In particular, the upper and lower shifting numbers coincide $\tau^+(\Phi)=\tau^-(\Phi)=\ell$.
\end{eg}

\begin{rmk}
\label{rmk:NotStandard}
Proposition~\ref{prop:EntropyBasic}\ref{propEntropy:conjugacyInvar} and Proposition~\ref{prop:ShiftingNumStrip} imply that the shifting numbers $\tau^\pm(\Phi)$ are conjugacy invariants. Therefore, if an autoequivalence $\Phi'$ has $\tau^+(\Phi')\neq\tau^-(\Phi')$, then it cannot be conjugate to any standard autoequivalence, by Example~\ref{eg:stdautoeq}. This will be used to show that spherical twists and $\bP$-twists are not conjugate to any standard autoequivalence.
\end{rmk}

$$
\sim\sim\sim
$$

Now, we recall a method for computing the shifting numbers, as developed in \cite{FanShifting}, using bounded $t$-structures on $\cD$. A full additive subcategory $\sA\subseteq\cD$ is the \emph{heart} of a bounded $t$-structure on $\cD$ if and only if:
\begin{enumerate}[label=(\roman*)]
\item $\Hom(A_1[k_1],A_2[k_2])=0$ if $k_1>k_2$ and $A_1,A_2\in\sA$,
\item for every nonzero object $E\in\cD$, there is a (unique) sequence of exact triangles
$$
\xymatrix@C=.5em{
 & 0 \ar[rrrr] &&&& E_1 \ar[rrrr] \ar[dll] &&&& E_2
\ar[rr] \ar[dll] && \cdots \ar[rr] && E_{\ell-1}
\ar[rrrr] &&&& E \ar[dll]  &   \\
&&& A_1[k_1] \ar@{-->}[ull] &&&& A_2[k_2] \ar@{-->}[ull] &&&&&&&& A_\ell[k_\ell] \ar@{-->}[ull] 
}
$$
with $k_1>\cdots>k_\ell$ and  $A_1,\ldots,A_\ell\in\sA\backslash\{0\}$.
\end{enumerate}
The object $A_i\in\sA\backslash\{0\}$ is called the cohomology object of $E$ at degree $k_i$, with respect to the heart $\sA$.

Throughout this article, the following notations will be used.

\begin{notn}
\label{notn:HeartCohomologyObjectPhase}
We denote the maximal (resp.~minimal) degrees and cohomology objects of $E$ with respect to the heart $\sA$ as follows:
$$
\phi^+_{\sA}(E)\coloneqq k_1, \quad \phi^-_{\sA}(E)\coloneqq k_\ell, \quad E_{\sA}^{+}\coloneqq A_1, \quad  E_{\sA}^{-}\coloneqq A_\ell.
$$
Additionally, the following ``cut-off" notations will also be used: For a real number $s$, suppose $k_p\geq s>k_{p+1}$, then we define
$$
E_\sA^{\geq s}\coloneqq E_p \qquad \text{ and } \qquad E_\sA^{<s}\coloneqq\Cone\left(E_p\rightarrow E\right).
$$
They satisfy $\phi_\sA^-(E_\sA^{\geq s})=k_p\geq s > k_{p+1}= \phi_\sA^+(E_\sA^{< s})$.
The notions $E_\sA^{>s}$ and $E_\sA^{\leq s}$ are defined similarly, and they satisfy $\phi_\sA^-(E_\sA^{> s})> s \geq \phi_\sA^+(E_\sA^{\leq s})$.
\end{notn}

The following proposition shows that the shifting numbers can be computed via the linear growth rates of the maximal (and minimal) degrees of cohomology objects of $\Phi^n(G)$ as $n\rightarrow\infty$.

\begin{prop}[\cite{FanShifting}*{Theorem~1.1}]
\label{prop:FanShiftingHeart}
Consider an endofunctor $\Phi\colon\cD\rightarrow\cD$ of a triangulated category $\cD$ with a split generator $G$. Let $\sA\subseteq\cD$ be the heart of a bounded $t$-structure.
\begin{enumerate}[label=(\roman*)]
    \item The limit 
    $$
    \lim_{n\rightarrow\infty}\frac{\phi^+_{\sA}(\Phi^nG)}{n}
    $$
    exists, is independent of the choices of $G$ and $\sA$, and coincides with $\tau^+(\Phi)$.
    \item The limit 
    $$
    \lim_{n\rightarrow\infty}\frac{\phi^-_{\sA}(\Phi^nG)}{n}
    $$
    exists, is independent of the choices of $G$ and $\sA$, and coincides with $\tau^-(\Phi)$.
\end{enumerate}
\end{prop}

$$
\sim\sim\sim
$$

Next, we recall the notion of \emph{categorical polynomial entropy function} $h_t^{\pol}(\Phi)$ \cite{FanFuOuchi}. It can be regarded as a more refined invariant: While the entropy function measures the \emph{exponential} growth rate of $\epsilon_t(G,\Phi^nG')$, the polynomial entropy function measures its \emph{polynomial} growth rate. It can be expressed as follows \cite{FanFuOuchi}*{Lemma~2.7}:
$$
h_t^\pol(\Phi)=\limsup_{n\rightarrow\infty}\frac{\log\epsilon_t(G,\Phi^n G')-nh_t(\Phi)}{\log(n)},
$$
and is independent of the choices of split generators $G,G'$.

$$
\sim\sim\sim
$$

Finally, we recall that when a triangulated category $\cD$ admits a \emph{stability condition} \cite{BriStab}, the complexity of an endofunctor can be measured by the growth rate of \emph{mass} with respect to the stability condition \cite{DHKK}*{Section~4.5}, \cite{Ikeda}. Let $\sigma=(Z,\sP)$ be a stability condition on $\cD$, where $Z\colon K_0(\cD)\rightarrow\bC$ is a group homomorphism, and $\sP=\{\sP(\phi)\}_{\phi\in\bR}$ is a collection of full additive subcategories of $\cD$, satisfying various axioms. For a nonzero object $E$, its \emph{mass function} with respect to $\sigma$ is the following real function in $t$:
$$
m_{\sigma,t}(E)=\sum_{k}|Z(A_k)|e^{\phi(A_k)t}
$$
where the $A_k$'s are the $\sigma$-semistable factors of $E$.
The \emph{mass growth function} of an endofunctor $\Phi\colon\cD\rightarrow\cD$ is defined as:
$$
h_{\sigma,t}(\Phi)=\limsup_{n\rightarrow\infty}\frac{1}{n}\log m_{\sigma,t}(\Phi^n(G)),
$$
and is independent of the choice of the split generator $G$ \cite{Ikeda}*{Theorem~3.5(1)}. 
Similarly, one can define the \emph{polynomial mass growth function} \cite{FanFuOuchi}*{Definition~3.3}:
$$
h_{\sigma,t}^\pol(\Phi)=\limsup_{n\rightarrow\infty}\frac{\log m_{\sigma,t}(\Phi^n(G))-nh_{\sigma,t}(\Phi)}{\log(n)}.
$$
The (polynomial) entropy functions and the (polynomial) mass growth functions are related as follows:
\begin{lemma}
\label{lemma:entropyVSmassgrowth}
Let $\Phi\colon\cD\rightarrow\cD$ be an endofunctor of a triangulated category. Then
\begin{enumerate}[label=(\roman*)]
    \item\label{item:entropyVSmassgrowth-1} $h_t(\Phi)\geq h_{\sigma,t}(\Phi)$. See \cite{Ikeda}*{Theorem~3.5(2)}.
    \item\label{item:entropyVSmassgrowth-2} For any $t\in\bR$, if $h_t(\Phi)=h_{\sigma,t}(\Phi)$, then $h_t^\pol(\Phi)\geq h_{\sigma,t}^\pol(\Phi)$. See \cite{FanFuOuchi}*{Lemma~3.6}.
\end{enumerate}
\end{lemma}

\begin{rmk}
\label{rmk:MassLowerBound}
We do not provide the full definition of stability conditions here. However, note that by the \emph{support property} of stability conditions, there exists a constant $C>0$ such that $|Z(M)|>C$ for any semistable object $M\in\sP_\sigma(\phi)$. This implies that if an object $E\in\cD$ has $\ell$ nonzero cohomology objects with respect to the heart $\sA_\sigma=\sP_\sigma(0,1]$, then we have a lower bound $m_{\sigma,0}(E)>\ell\cdot C$.
\end{rmk}

\subsection{Spherical twists}
In this subsection, we recall the notion of \emph{spherical twists} \cite{SeidelThomas}, and previous results on the categorical entropy functions of spherical twists \cite{OuchiSpherical}.

An object $E\in\cD$ is called a \emph{$d$-spherical object} if $\bS(E)\cong E[d]$ and
$$
\dim\Hom(E,E[k])=
\begin{cases}
    1 & \text{ if } k=0 \text{ or }d, \\
    0 & \text{ otherwise}.
\end{cases}
$$
Examples of spherical objects include line bundles in the derived categories of Calabi--Yau varieties and structure sheaves of rational $(-2)$-curves in projective surfaces, among others.

One of the most interesting features of spherical objects is that they induce autoequivalences of the triangulated category, known as \emph{spherical twists} \cite{SeidelThomas}. For a $d$-spherical object $E$, the associated spherical twist, denoted by $T_E\colon\cD\rightarrow\cD$, sends an object $F$ to
$$
T_E(F)=\Cone\left(\RHom(E,F)\otimes E \xrightarrow{\text{ev}}F\right).
$$
It follows that $T_E(E)\cong E[1-d]$, and $T_E(F)\cong F$ for any $F\in E^\perp=\left\{M\in\cD:\RHom(E,M)=0\right\}$.

\begin{thm}[\cite{OuchiSpherical}*{Theorem~3.1}]
\label{thm:Ouchispherical}
Let $E\in\cD$ be a $d$-spherical object, where $d\geq2$. The categorical entropy function of the spherical twist $T_E$ is given by:
\begin{enumerate}[label=(\roman*)]
    \item\label{item:Genki-1} $h_t(T_E)=(1-d)t$ for $t\leq0$.
    \item\label{item:Genki-2} $h_t(T_E)\leq0$ for $t>0$.
    \item If $E^\perp\neq\{0\}$, then $h_t(T_E)=0$ for $t>0$.
\end{enumerate}
\end{thm}

We refer to \cite{OuchiSpherical} for the proof.
In Section~\ref{subsubsec:spherical}, we will show that the condition $E^\perp\neq\{0\}$ can be replaced by certain conditions that are always satisfied when $\cD=\Db(X)$. The following lemma provides an example of such weaker conditions.

\begin{lemma}
\label{lemma:SphericalMaxDegree}
Let $E$ be a $d$-spherical object, where $d\geq2$. Suppose there exists a nonzero object $F\in\cD$ and a heart $\sA\subseteq\cD$ of bounded $t$-structure such that
$$
\lim_{n\rightarrow\infty}\frac{\phi_{\sA}^+(T_E^n(F))}{n}=0.
$$
Then $h_t(T_E)=0$ for $t>0$.
\end{lemma}
\begin{proof}
Let $G$ be a split generator of $\cD$. Then $G\oplus F$ is also a split generator. By Proposition~\ref{prop:FanShiftingHeart}, 
$$
\tau^+(T_E)  = \lim_{n\rightarrow\infty} \frac{\phi_\sA^+(T_E^n(G\oplus F))}{n}
\geq \lim_{n\rightarrow\infty} \frac{\phi_\sA^+(T_E^n(F))}{n} =0.
$$
From Theorem~\ref{thm:Ouchispherical}\ref{item:Genki-1}, we know that $h_0(T_E)=0$. Thus, for $t>0$, Proposition~\ref{prop:ShiftingNumStrip} gives
$$
h_t(T_E)=t\cdot \tau^+(T_E)\geq0.
$$
Combining this with theorem~\ref{thm:Ouchispherical}\ref{item:Genki-2}, we conclude that $h_t(T_E)=0$ for all $t>0$.
\end{proof}

\begin{rmk}
Note that if the condition $E^\perp\neq\{0\}$ holds, then any $0\neq F\in E^\perp$ will satisfy $\lim_{n\rightarrow\infty}\frac{\phi_{\sA}^+(T_E^n(F))}{n}=0$, since  $T_E^n(F)=F$ for all $n$. In Section~\ref{subsubsec:spherical}, we will establish conditions that ensure the existence of such an object $F$ and heart $\sA$ which, while not necessarily meeting the condition as strong as $T_E^n(F)=F$, still satisfies $\lim_{n\rightarrow\infty}\frac{\phi_{\sA}^+(T_E^n(F))}{n}=0$.
\end{rmk}

\subsection{$\bP$-twists}
\label{subsec:P-twistOuchiVersion}
In this subsection, we recall the notion of \emph{$\bP$-twists} \cite{HuyThoPtwist}. We compute their categorical entropy functions, and obtain results analogous to Theorem~\ref{thm:Ouchispherical} and Lemma~\ref{lemma:SphericalMaxDegree}.

An object $E\in\cD$ is called a \emph{$\bP^d$-object} if $\bS(E)\cong E[2d]$ and $\Hom(E,E[\ast])\cong H^*(\bC\bP^d,\bZ)\otimes\mathbf{k}$ as $\mathbf{k}$-algebras. Examples of $\bP$-objects include line bundles in the derived categories of hyperk\"ahler manifolds and the structure sheaf of an embedded $\bP^d$ in a $2d$-dimensional holomorphic symplectic variety, among others.

Similar to spherical objects, a $\bP$-object also induces an autoequivalence of the triangulated category, known as the \emph{$\bP$-twist}. We now recall its definition following \cite{HuyThoPtwist}. For a $\bP$-object $E$, a generator $h\in\Hom(E,E[2])$ can be viewed as a morphism $h\colon E[-2]\rightarrow E$. Denote the image of $h$ under the natural isomorphism $\Hom(E,E[2])\cong\Hom(E^\vee,E^\vee[2])$ by $h^\vee$. The associated $\bP$-twist $P_E\colon\cD\rightarrow\cD$ sends an object $F$ to
$$
P_E(F)=\Cone\left(\Cone\left(\Hom^{\ast-2}(E,F)\otimes E\xrightarrow{h^\vee\cdot\id-\id\cdot h}\Hom^\ast(E,F)\otimes E\right)\rightarrow F\right).
$$
It follows that $P_E(E)\cong E[-2d]$, and $P_E(F)\cong F$ for any $F\in E^\perp$ (see \cite{HuyThoPtwist}*{Lemma~2.5}).

\begin{thm}
\label{thm:OuchiP-twistVersion}
Let $E\in\cD$ be a $\bP^d$-object. The categorical entropy function of the $\bP$-twist $P_E$ is given by:
\begin{enumerate}[label=(\roman*)]
    \item $h_t(P_E)=-2dt$ for $t\leq0$.
    \item $h_t(P_E)\leq0$ for $t>0$.
    \item If $E^\perp\neq\{0\}$, then $h_t(P_E)=0$ for $t>0$.
\end{enumerate}
\end{thm}

\begin{proof}
Let $G$ be a split generator of $\cD$. Define
$$
M\coloneqq\Cone\left(\Hom^{\ast-2}(E,G)\otimes E\xrightarrow{h^\vee\cdot\id-\id\cdot h}\Hom^\ast(E,G)\otimes E\right)[1].
$$
Applying $P_E^{n-1}$ to the exact triangle
$$
G\rightarrow P_E(G)\rightarrow M \xrightarrow{+1},
$$
one obtains
$$
 P_E^{n-1}(G)\rightarrow P_E^n(G)\rightarrow M[-2d(n-1)]\xrightarrow{+1}.
$$
Then
\begin{align*}
\epsilon_t(G,P_E^n(G)) & \leq \epsilon_t(G,M) e^{-2d(n-1)t} + \epsilon_t(G,P_E^{n-1}(G)) \\
& \leq \epsilon_t(G,M) e^{-2d(n-1)t} + \epsilon_t(G,M) e^{-2d(n-2)t} + \epsilon_t(G,P_E^{n-2}(G)) \\
& \leq \cdots \\
& \leq \epsilon_t(G,M)\left(\sum_{k=0}^{n-1} e^{-2dkt}\right) +\epsilon_t(G,G).
\end{align*}
First, let us consider the case for $t\leq0$.
\begin{align*}
h_t(P_E) & =\lim_{n\rightarrow\infty}\frac{1}{n}\log\epsilon_t(G,P_E^n(G)) \\
& \leq \lim_{n\rightarrow\infty}\frac{1}{n}\log \left(\epsilon_t(G,M)\left(\sum_{k=0}^{n-1} e^{-2dkt}\right) +\epsilon_t(G,G)\right)\\
& \leq \lim_{n\rightarrow\infty}\frac{1}{n}\log \left(\epsilon_t(G,M)\cdot n \cdot e^{-2d(n-1)t}+\epsilon_t(G,G)\right)\leq -2dt.
\end{align*}
On the other hand, since $G\oplus E$ is a split generator of $\cD$, 
\begin{align*}
h_t(P_E) & =\lim_{n\rightarrow\infty}\frac{1}{n}\log\epsilon_t(G\oplus E,P_E^n(G\oplus E)) \\
& \geq \lim_{n\rightarrow\infty}\frac{1}{n}\log\epsilon_t(G\oplus E,P_E^n(E)) \\
& = \lim_{n\rightarrow\infty}\frac{1}{n}\log \left(\epsilon_t(G\oplus E,E)e^{-2dnt}\right)=-2dt.
\end{align*}
This proves $h_t(P_E)=-2dt$ for $t\leq0$.

Next, we consider the case for $t>0$.
\begin{align*}
h_t(P_E) & =\lim_{n\rightarrow\infty}\frac{1}{n}\log\epsilon_t(G,P_E^n(G)) \\
& \leq \lim_{n\rightarrow\infty}\frac{1}{n}\log \left(\epsilon_t(G,M)\left(\sum_{k=0}^{n-1} e^{-2dkt}\right) +\epsilon_t(G,G)\right)\\
& \leq \lim_{n\rightarrow\infty}\frac{1}{n}\log \left(\epsilon_t(G,M)\cdot n+\epsilon_t(G,G)\right)\leq 0.
\end{align*}
Moreover, assuming that $E^\perp\neq\{0\}$, say $0\neq F\in E^\perp$. Consider the split generator $G\oplus F$,
\begin{align*}
h_t(P_E) & =\lim_{n\rightarrow\infty}\frac{1}{n}\log\epsilon_t(G\oplus F,P_E^n(G\oplus F)) \\
& \geq \lim_{n\rightarrow\infty}\frac{1}{n}\log\epsilon_t(G\oplus F,P_E^n(F)) \\
& = \lim_{n\rightarrow\infty}\frac{1}{n}\log \epsilon_t(G\oplus F,F)=0.
\end{align*}
\end{proof}

In Section~\ref{subsubsec:P-twist}, we will show that, similar to spherical twists, the condition $E^\perp\neq\{0\}$ can be replaced by certain conditions that are always satisfied when $\cD=\Db(X)$.
The following lemma is the $\bP$-twist analog of Lemma~\ref{lemma:SphericalMaxDegree}, and it can be proved using the same argument.

\begin{lemma}
\label{lemma:P-MaxDegree}
Let $E$ be a $\bP$-object. Suppose there exists a nonzero object $F\in\cD$ and a heart $\sA\subseteq\cD$ of bounded $t$-structure such that
$$
\lim_{n\rightarrow\infty}\frac{\phi_{\sA}^+(P_E^n(F))}{n}=0.
$$
Then $h_t(P_E)=0$ for $t>0$.
\end{lemma}

\section{Proof of Main Theorems}
In Section~\ref{subsec:Rigid(cde)}, we introduce various conditions that will be imposed for computing the categorical (polynomial) entropy functions of spherical twists and $\bP$-twists.
In Section~\ref{subsec:ProofOfThm1}, we study the categorical entropy functions and prove Theorem~\ref{thm:MainP-twist}.
In Section~\ref{subsec:ProofOfThm2}, we study the categorical polynomial entropy functions and prove Theorem~\ref{thm:PolyCat}.

\subsection{Rigid objects and Conditions (c), (d), (e)}
\label{subsec:Rigid(cde)}

The conditions in the following definition are satisfied by both spherical and $\bP$-objects.

\begin{defn}
Let $d\geq2$ be a positive integer. An object $E\in\cD$ is called \emph{$d$-rigid} if it satisfies the following conditions:
\begin{enumerate}[label=(\alph*)]
\item\label{cond:Serre} $\bS(E)\cong E[d]$,
\item\label{cond:Hom(E)} $\Hom(E,E)=\bC$, $\Hom(E,E[1])=0$, and $\Hom(E,E[k])=0$ for $k<0$.
\setcounter{nameOfYourChoice}{\value{enumi}}
\end{enumerate}
\end{defn}

\noindent Note that a $d$-rigid object $E$ also has the vanishing $\Hom(E,E[k])=0$  for $k>d$, since
$$
\Hom(E,E[k])\cong\Hom(E[k],\bS(E))^\vee \cong \Hom(E[k],E[d])^\vee \cong \Hom(E,E[d-k])^\vee.
$$
In the following, we gather some conditions on $d$-rigid objects that will be imposed at various stages of the computation of entropy functions.

\begin{defn}
Let $E\in\cD$ be a $d$-rigid object. 
We say $E$ satisfies \emph{Condition~\ref{cond:hasHomN0}} if:
\begin{enumerate}[label=(\alph*)]
    \setcounter{enumi}{\value{nameOfYourChoice}}
    \item\label{cond:hasHomN0} $\cD$ is indecomposable and $\left<E\right>\subsetneq\cD$. (Here, $\left<E\right>$ denotes the smallest triangulated subcategory of $\cD$ containing $E$.)
    \setcounter{nameOfYourChoice}{\value{enumi}}
\end{enumerate}
We say $E$ satisfies \emph{Condition~\ref{cond:existHeartPerp}} if:
\begin{enumerate}[label=(\alph*)]
    \setcounter{enumi}{\value{nameOfYourChoice}}
    \item\label{cond:existHeartPerp}  there exists a heart $\sA\subseteq\cD$ of bounded $t$-structure such that (at least) one of the following holds:
        \begin{itemize}
            \item there exists a nonzero object $A\in\sA$ such that $\Hom(E_{\sA}^{+},A)=0$, or
            \item $\phi_\sA^+(E)-\phi_\sA^-(E)\leq1$.
        \end{itemize}
    (The notations are defined in Notation~\ref{notn:HeartCohomologyObjectPhase}.)
    \setcounter{nameOfYourChoice}{\value{enumi}}
\end{enumerate}
We say $E$ satisfies \emph{Condition~\ref{cond:StabilityHeart}} if:
\begin{enumerate}[label=(\alph*)]
    \setcounter{enumi}{\value{nameOfYourChoice}}
    \item\label{cond:StabilityHeart} there exists a Bridgeland stability condition $\sigma=(Z_\sigma,\sP_\sigma)$ on $\cD$ such that $E\in\sP_\sigma(-1,1]$.
    \setcounter{nameOfYourChoice}{\value{enumi}}
\end{enumerate}
\end{defn}

\begin{rmk}
\label{rmk:(c)(d)D(X)satisfy}
Note that Condition~\ref{cond:hasHomN0} is satisfied for $\cD=\Db(X)$ when $X$ is connected, which we always assume.
Condition~\ref{cond:existHeartPerp} is also satisfied for $\cD=\Db(X)$: In this case, we can simply take $\sA=\Coh(X)$ and use the fact that for any coherent sheaf $E_{\sA}^{+}\in\sA$, we have $\Hom(E_{\sA}^{+},\cO(-n))=0$ for sufficiently large $n$.

Condition~\ref{cond:StabilityHeart} is more restrictive than Condition~\ref{cond:existHeartPerp}, and will be imposed only in the computation of categorical polynomial entropy $h_0^\pol(-)$ of spherical twists and $\bP$-twists.
\end{rmk}

\subsection{Entropy functions and Proof of Theorem~\ref{thm:MainP-twist}}
\label{subsec:ProofOfThm1}

\subsubsection{Lemmas}

We begin with an elementary lemma concerning the extremal degrees $\phi_\sA^\pm(-)$ and cohomology objects $E_\sA^\pm$ of objects in an exact triangle. The relevant notations are defined in Notation~\ref{notn:HeartCohomologyObjectPhase}.

\begin{lemma}
\label{lemma:ExtremalDegree}
Let $\sA\subseteq\cD$ be the heart of a bounded $t$-structure on a triangulated category $\cD$.
Suppose $M\rightarrow E\rightarrow N\xrightarrow{+1}$ is an exact triangle, where $E,M,N$ are nonzero objects in $\cD$.
Then:
\begin{enumerate}[label=(\roman*)]
    \item\label{item:ExtremalDegree-1}$\min\{\phi_{\sA}^-(M),\phi_{\sA}^-(N)\}\leq\phi_{\sA}^-(E)\leq\phi_{\sA}^+(E)\leq\max\{\phi_{\sA}^+(M),\phi_{\sA}^+(N)\}$.
    \item\label{item:NonzeroCohom} If $\phi_\sA^-(M)\geq\phi_\sA^+(N)$, then $E$ has nonzero cohomology objects at degrees $\phi_\sA^-(M)$ and $\phi_\sA^+(N)$.
    \item\label{item:ExtremalDegree-2} If $\phi_{\sA}^+(M)\geq\phi_{\sA}^+(N)$, then $\phi_{\sA}^+(E)=\phi_{\sA}^+(M)$. Similarly, if $\phi_{\sA}^-(N)\leq\phi_{\sA}^-(M)$, then $\phi_{\sA}^-(E)=\phi_{\sA}^-(N)$.
    \item\label{item:ExtremalDegree-3} If either of the following holds:
        \begin{itemize}
            \item $\phi_{\sA}^+(N)\geq\phi_{\sA}^+(M)+2$, or
            \item $\phi_{\sA}^+(N)=\phi_{\sA}^+(M)+1$ and $\Hom(N_{\sA}^{+},M_{\sA}^{+})=0$,
        \end{itemize}
    then $\phi_{\sA}^+(E)=\phi_{\sA}^+(N)$ and $E_\sA^+\cong N_\sA^+$.
    \item\label{item:ExtremalDegree-4} Similarly, if either of the following holds:
        \begin{itemize}
            \item $\phi_{\sA}^-(M)\leq\phi_{\sA}^-(N)-2$, or
            \item $\phi_{\sA}^-(M)=\phi_{\sA}^-(N)-1$ and $\Hom(N_{\sA}^{-},M_{\sA}^{-})=0$,
        \end{itemize}
    then $\phi_{\sA}^-(E)=\phi_{\sA}^-(M)$ and $E_\sA^-\cong M_\sA^-$.
\end{enumerate}
\end{lemma}
\begin{proof}[Proof of \ref{item:ExtremalDegree-1}]\renewcommand{\qedsymbol}{}
This follows directly from that 
$$
\Hom\left(E_{\sA}^{+}[\phi_{\sA}^+(E)],E\right)\neq0 \quad \text{ and } \quad \Hom\left(E,E_{\sA}^{-}[\phi_{\sA}^-(E)]\right)\neq0.
$$
\end{proof}

\begin{proof}[Proof of \ref{item:NonzeroCohom}]\renewcommand{\qedsymbol}{}
First, if $\phi_{\sA}^-(M)>\phi_{\sA}^+(N)$, then the cohomology filtrations of $M$ and $N$ (with respect to $\sA$) can be combined to form the cohomology filtration of $E$:

\resizebox{15cm}{!}{
\xymatrix@C=.5em{
0 \ar[rrrr]
&&&&
\ast \ar[rr] \ar[dll]
&&
\cdots \ar[rr]
&&
\ast \ar[rrrr]
&&&&
\ast \ar[rrrr] \ar[dll]
&&&&
\ast \ar[rr] \ar[dll]
&&
\cdots \ar[rr]
&&
\ast \ar[rrrr]
&&&&
E \ar[dll]
\\
&&
M_\sA^+[\phi_\sA^+(M)] \ar@{-->}[ull]
&&&&&&&&
M_\sA^-[\phi_\sA^-(M)] \ar@{-->}[ull]
&&&&
N_\sA^+[\phi_\sA^+(N)] \ar@{-->}[ull]
&&&&&&&&
N_\sA^-[\phi_\sA^-(N)] \ar@{-->}[ull]
}
}
\vspace{.5cm}

\noindent Therefore, $E$ has nonzero cohomology objects at both degrees $\phi_\sA^-(M)$ and $\phi_\sA^+(N)$.

Second, suppose $\phi_{\sA}^-(M)=\phi_{\sA}^+(N)=\phi$. Then we have
$$
\xymatrix@C=.5em{
0 \ar[rrrr]
&&&&
\ast \ar[rrrr] \ar[dll]
&&&&
\ast \ar[rrrr] \ar[dll]
&&&&
\ast \ar[rrrr] \ar[dll]
&&&&
E \ar[dll]
\\
&&
M_\sA^{>\phi} \ar@{-->}[ull]
&&&&
M_\sA^-[\phi] \ar@{-->}[ull]
&&&&
N_\sA^+[\phi] \ar@{-->}[ull]
&&&&
N_\sA^{<\phi} \ar@{-->}[ull]
}
$$
Thus, $E$ has a nonzero cohomology object at degree $\phi$, which is an extension of $M_\sA^-$ and $N_\sA^+$ in $\sA$.
\end{proof}

\begin{proof}[Proof of \ref{item:ExtremalDegree-2}]\renewcommand{\qedsymbol}{}
The object $M$ sits in the exact triangle
$$
M_{\sA}^+[\phi_\sA^+(M)] \rightarrow M \rightarrow M_{\sA}^{<\phi_\sA^+(M)} \xrightarrow{+1}.
$$
Thus, there  exists an object $X$ and exact triangles
$$
M_{\sA}^{<\phi_\sA^+(M)}\rightarrow X\rightarrow N \xrightarrow{+1} \quad \text{ and } \quad M_{\sA}^+[\phi_\sA^+(M)] \rightarrow E \rightarrow X \xrightarrow{+1}.
$$
Applying \ref{item:ExtremalDegree-1} to the first exact triangle, we have 
$$
\phi_\sA^+(X)\leq\max\left\{\phi_{\sA}^+\left(M_{\sA}^{<\phi_\sA^+(M)}\right),\phi_\sA^+(N)\right\} \leq \phi_\sA^+(M).
$$
Then, applying \ref{item:NonzeroCohom} to the second exact triangle, we obtain $\phi_\sA^+(E)=\phi_\sA^+(M)$.
The statement concerning $\phi_\sA^-(-)$ can be proved similarly.
\end{proof}

\begin{proof}[Proof of \ref{item:ExtremalDegree-3}]
There are exact triangles
$$
\xymatrix@C=.5em{
M \ar[rrrr]
&&&&
X \ar[rrrr] \ar[dll]
&&&&
E \ar[dll]
\\
&&
N_\sA^+[\phi_\sA^+(N)] \ar@{-->}[ull]
&&&&
N_\sA^{<\phi_\sA^+(N)} \ar@{-->}[ull]
}
$$
First, suppose $\phi_{\sA}^+(N)\geq\phi_{\sA}^+(M)+2$. Then $\Hom\left(N_\sA^+[\phi_\sA^+(N)],M[1]\right)=0$.
Therefore, $X\cong M\oplus N_\sA^+[\phi_\sA^+(N)]$.
Consequently, the order of $M$ and $N_\sA^+[\phi_\sA^+(N)]$ can be exchanged, resulting in $\phi_\sA^+(E)=\phi_\sA^+(N)$ and $E_\sA^+=N_\sA^+$.

Second, suppose $\phi_{\sA}^+(N)=\phi_{\sA}^+(M)+1$ and $\Hom(N_{\sA}^{+},M_{\sA}^{+})=0$. There are exact triangles
$$
\xymatrix@C=.5em{
0 \ar[rrrr]
&&&&
\ast \ar[rrrr] \ar[dll]
&&&&
\ast \ar[rrrr] \ar[dll]
&&&&
\ast \ar[rrrr] \ar[dll]
&&&&
E \ar[dll]
\\
&&
M_\sA^+[\phi_\sA^+(M)] \ar@{-->}[ull]
&&&&
M_\sA^{<\phi_\sA^+(M)} \ar@{-->}[ull]
&&&&
N_\sA^+[\phi_\sA^+(N)] \ar@{-->}[ull]
&&&&
N_\sA^{<\phi_\sA^+(N)} \ar@{-->}[ull]
}
$$
Consider the second and third objects from the left.
Since $\Hom\left(N_\sA^+[\phi_\sA^+(N)],M_\sA^{<\phi_\sA^+(M)}[1]\right)=0$, we can exchange their order, moving $N_\sA^+[\phi_\sA^+(N)]$ to the second position.
Next, we have $\Hom\left(N_\sA^+[\phi_\sA^+(N)],M_\sA^+[\phi_\sA^+(M)+1]\right)=0$ by assumption. Therefore, we can also exchange their order, moving $N_\sA^+[\phi_\sA^+(N)]$ to the first position, and obtain
$$
\xymatrix@C=.5em{
0 \ar[rrrr]
&&&&
\ast \ar[rrrr] \ar[dll]
&&&&
\ast \ar[rrrr] \ar[dll]
&&&&
\ast \ar[rrrr] \ar[dll]
&&&&
E \ar[dll]
\\
&&
N_\sA^+[\phi_\sA^+(N)] \ar@{-->}[ull]
&&&&
M_\sA^+[\phi_\sA^+(M)] \ar@{-->}[ull]
&&&&
M_\sA^{<\phi_\sA^+(M)} \ar@{-->}[ull]
&&&&
N_\sA^{<\phi_\sA^+(N)} \ar@{-->}[ull]
}
$$
Consequently, we have $\phi_\sA^+(E)=\phi_\sA^+(N)$ and $E_\sA^+=N_\sA^+$.

Finally, \ref{item:ExtremalDegree-4} can be proved similarly.
\end{proof}

The following lemma will be useful for managing the degrees of nonzero cohomology objects.

\begin{lemma}
\label{lemma:concentrateInd-1}
Let $E\in\cD$ be a $d$-rigid object satisfying Condition~\ref{cond:hasHomN0} for some $d\geq2$.
There exists $F\in\cD$ such that:
    \begin{itemize}
        \item $\Hom(E,F)\neq0$, and
        \item $\Hom(E,F[k])=0$ if $k<0$ or $k>d-2$.
    \end{itemize}
\end{lemma}

\begin{proof}
First, note that Condition~\ref{cond:hasHomN0} implies that $\left<E,E^\perp\right>\subsetneq\cD$. Indeed, assume the contrary, that $\left<E,E^\perp\right>=\cD$. There are two possibilities:
    \begin{enumerate}[label=(\roman*)]
        \item If $E^\perp=\{0\}$, then $\cD$ is generated by $E$.
        \item If $E^\perp\neq\{0\}$, then $\cD\cong\left<E\right>\oplus E^\perp$, since for any $A\in E^\perp$, we have
        $$
        \Hom(A,E[k])\cong\Hom(A,\bS(E)[k-d])\cong\Hom(E[k-d],A)^\vee=0 \quad \text{ for all } k.
        $$
    \end{enumerate}

Let $F_1\in\cD$ be an object such that $F_1\notin\left<E,E^\perp\right>$.
By shifting, we may assume that $\Hom(E,F_1)\neq0$ and $\Hom(E,F_1[k])=0$ for all $k<0$. Denote
$$
n=\max\left\{k:\Hom(E,F_1[k])\neq0\right\}.
$$
If $n\leq d-2$, then we are done. Now suppose that $n\geq d-1$. Let
$$
F_2=\Cone\left(\Hom(E,F_1)\otimes E\xrightarrow{\text{ev}}F_1\right)[1].
$$
Applying $\Hom(E,-)$ to the exact triangle
$$
\Hom(E,F_1)\otimes E\xrightarrow{\text{ev}}F_1\rightarrow F_2[-1]\xrightarrow{+1},
$$
one obtains
\begin{center}
\begin{tabular}{ c c c c c c }
$\rightarrow$ & $0$ & $\rightarrow$ & $0$ & $\rightarrow$ & $\Hom(E,F_2[-2])$ \\
$\rightarrow$ & $\Hom(E,F_1)\otimes\Hom(E,E)$ & $\xrightarrow{\cong}$ & $\Hom(E,F_1)$ & $\rightarrow$ & $\Hom(E.F_2[-1])$ \\
$\rightarrow$ & $0$ & $\rightarrow$ & $\Hom(E,F_1[1])$ & $\rightarrow$ & $\Hom(E,F_2)$ \\
$\rightarrow$ & $\Hom(E,F_1)\otimes\Hom(E,E[2])$ & $\rightarrow$ & $\Hom(E,F_1[2])$ & $\rightarrow$ & $\Hom(E.F_2[1])$ \\
$\rightarrow$ & $\cdots$
\end{tabular}
\end{center}
where the isomorphism in the second row and the zero in the third row both follow from Condition~\ref{cond:Hom(E)} of $E$ being a $d$-rigid object.
Therefore, $\Hom(E,F_2[k])=0$ for all $k<0$ and $k>n-1$.
Note that $F_2\notin\left<E,E^\perp\right>$; otherwise, we would have $F_1\in\left<E,E^\perp\right>$, which is a contradiction.
Thus, we can find the desired object $F$ inductively.
\end{proof}

\subsubsection{Entropy function of spherical twists}
\label{subsubsec:spherical}

\begin{prop}
\label{prop:SphericalMaxDegree}
Let $E$ be a $d$-spherical object in $\cD$, where $d\geq2$.
Suppose $E$ satisfies Conditions~\ref{cond:hasHomN0} and \ref{cond:existHeartPerp}. Then, there exists a nonzero object $F\in\cD$ and a heart $\sA\subseteq\cD$ such that
$$
\lim_{n\rightarrow\infty}\frac{\phi_{\sA}^+(T_E^n(F))}{n}=0.
$$
In fact, there exists a nonzero object $\widetilde{F}$ and an integer $\ell>0$ such that for any $n>\ell$, there is an exact triangle
$$
\widetilde{F}\rightarrow T_E^n(F) \rightarrow C_n \xrightarrow{+1}
$$
with $\phi_{\sA}^-(\widetilde{F})>\phi_{\sA}^+(C_n)$.
Therefore, $\phi_\sA^+(T_E^n(F))$ stabilizes when $n>\ell$.
\end{prop}

\begin{proof}
Let $F$ be a nonzero object in $\cD$. Applying $T_E^{n-1}$ to the exact triangle
$$
F\rightarrow T_E(F)\rightarrow \RHom(E,F)\otimes E[1] \xrightarrow{+1},
$$
one obtains
$$
T_E^{n-1}(F)\rightarrow T_E^n(F) \rightarrow \RHom(E,F)\otimes E[1+(n-1)(1-d)] \xrightarrow{+1}.
$$
Denote $M\coloneqq\RHom(E,F)\otimes E[1]$. Then we have
$$
\xymatrix@C=.5em{
F \ar[rrrr]
&&&&
\ast \ar[rrrr] \ar[dll]
&&&&
\ast \ar[rr] \ar[dll]
&&
\cdots \ar[rr] 
&&
\ast \ar[rrrr]
&&&&
T_E^n(F) \ar[dll] 
\\
&&
M \ar@{-->}[ull]
&&&&
M[1-d] \ar@{-->}[ull]
&&&&&&&&
M[(n-1)(1-d)] \ar@{-->}[ull]
}
$$
Since $E$ satisfies Condition~\ref{cond:existHeartPerp}, there exists a heart $\sA\subseteq\cD$ such that (at least) one of the following statements holds:
\begin{itemize}
    \item there exists a nonzero object $A\in\sA$ such that $\Hom(E_{\sA}^{+},A)=0$, or
    \item $\phi_\sA^+(E)-\phi_\sA^-(E)\leq1$.
\end{itemize}

First, suppose there exists $A\in\sA$ such that $\Hom(E_{\sA}^{+},A)=0$.
We claim that taking $F=A$ would satisfy the desired properties stated in the proposition. If $\RHom(E,F)=0$, then $T_E^n(F)=F$ for all $n$, and we can simply take $\widetilde{F}=F$ and $C_n=0$.
Otherwise, we have $M\neq0$ and $M_{\sA}^{+}\cong (E_{\sA}^{+})^{\oplus p}$ for some integer $p\geq1$.
Now, we compare $\phi_\sA^+(F)(=\phi_\sA^-(F)=0)$ and $\phi_\sA^+(M)$.
\begin{enumerate}[label=(\roman*)]
    \item If $\phi_\sA^+(M)\leq0$, then $\phi_\sA^+(T_E^n(F))=0$ by Lemma~\ref{lemma:ExtremalDegree}\ref{item:ExtremalDegree-2}. Moreover, $(T_E^n(F))_\sA^+$ is independent of $n$, since $\phi_\sA^+(M[k(1-d)])<0$ for any $k>0$. Therefore, taking $\widetilde{F}=(T_E(F))_\sA^+$ and $\ell=1$ would satisfy the desired property.
    \item If $\phi_{\sA}^+(M)\geq2$, then $\phi_{\sA}^+(T_E^n(F))=\phi_{\sA}^+(M)$ by Lemma~\ref{lemma:ExtremalDegree}\ref{item:ExtremalDegree-3}. Moreover, $(T_E^n(F))_\sA^+\cong M_\sA^+$ for all $n$. Therefore, taking $\widetilde{F}=M_\sA^+$ and $\ell=1$ would satisfy the desired property.
    \item If $\phi_{\sA}^+(M)=1$, using $\Hom(M_{\sA}^{+},F)\cong\Hom\left((E_{\sA}^{+})^{\oplus p},A\right)=0$, we still have $\phi_{\sA}^+(T_E^n(F))=\phi_{\sA}^+(M)$ and $(T_E^n(F))_\sA^+\cong M_\sA^+$ for all $n$, by Lemma~\ref{lemma:ExtremalDegree}\ref{item:ExtremalDegree-3}.
\end{enumerate}

Second, suppose $\phi_\sA^+(E)-\phi_\sA^-(E)\leq1$.
By Lemma~\ref{lemma:concentrateInd-1}, there exists $F\in\cD$ such that $\Hom(E,F)\neq0$, and $\Hom(E,F[k])=0$ if $k<0$ or $k>d-2$. Therefore, $M=\RHom(E,F)\otimes E[1]$ is nonzero and satisfies
$$
\phi_\sA^+(M)-\phi_\sA^-(M)\leq (d-2)+1=d-1.
$$
Choose $\ell>0$ large enough so that
$$
\phi_\sA^-(F)>\phi_\sA^+(M[\ell(1-d)])=\phi_\sA^+(M)+\ell(1-d).
$$
We claim that $T_E^n(F)$ admits a nonzero cohomology object at degree $\phi_\sA^+(M)+\ell(1-d)$ for any $n>\ell$.
Moreover, $\left(T_E^n(F)\right)_\sA^{\geq \phi_\sA^+(M)+\ell(1-d)}$ is independent of $n>\ell$.
To visualize this,  let
$$
\widetilde{\phi}\coloneqq\phi_\sA^+(M)+\ell(1-d),
$$
and consider the following picture:
\begin{center}
\begin{tikzpicture}
\draw (-1.5, 0) rectangle  ++(3,1) node[pos=.5, align=center] {$\cdots$};
\draw (1.5, 0) rectangle  ++(1.5,1) node[pos=.5, align=center] {$\sA[\widetilde{\phi}]$};
\draw (3, 0) rectangle  ++(2,1) node[pos=.5, align=center] {$\cdots$};
\draw (5, 0) rectangle  ++(3,1) node[pos=.5, align=center] {$\sA[\widetilde{\phi}+(1-d)]$};
\draw (8, 0) rectangle  ++(2,1) node[pos=.5, align=center] {$\cdots$};
\draw (10, 0) rectangle  ++(3,1) node[pos=.5, align=center] {$\sA[\widetilde{\phi}+2(1-d)]$};
\draw (13, 0) rectangle  ++(2,1) node[pos=.5, align=center] {$\cdots$};
\filldraw[draw=black,fill=yellow] (-1.5,-1)rectangle ++(2.5,1) node[pos=.5, align=center] {$F$};
\filldraw[draw=black,fill=yellow] (-1.5,-2)rectangle ++(4.5,1) node[pos=.5, align=center] {$M[(\ell-1)(1-d)]$};
\filldraw[draw=black,fill=yellow] (1.5,-3)rectangle ++(6.5,1) node[pos=.5, align=center] {$M[\ell(1-d)]$};
\draw (2.25,-2.5) node {$\times$};
\filldraw[draw=black,fill=yellow] (5,-4)rectangle ++(8,1) node[pos=.5, align=center] {$M[(\ell+1)(1-d)]$};
\draw (6.5,-3.5) node {$\times$};
\filldraw[draw=black,fill=yellow] (10,-5)rectangle ++(5,1);
\draw (11.5,-4.5) node {$\times$};
\end{tikzpicture}
\end{center}
The yellow regions indicate the possible ranges of nonzero cohomological degrees (with respect to the heart $\sA$) of the objects. The symbol ``$\times$" indicates that the cohomology object at that particular degree must be nonzero.

Observe that
$$
\min\left\{\phi_\sA^-(F),\phi_\sA^-(M),\phi_\sA^-(M[1-d]),\ldots,\phi_\sA^-(M[(\ell-1)(1-d)])\right\}\geq\widetilde{\phi}=\phi_\sA^+(M[\ell(1-d)]),
$$
and
$$
\max\left\{M[(\ell+1)(1-d)],M[(\ell+2)(1-d)],\ldots\right\}<\widetilde{\phi}.
$$
Thus, by Lemma~\ref{lemma:ExtremalDegree}\ref{item:NonzeroCohom}, $T_E^n(F)$ admits a nonzero cohomology object at degree $\widetilde{\phi}=\phi_\sA^+(M)+\ell(1-d)$ for any $n>\ell$.
Moreover, $\left(T_E^n(F)\right)_\sA^{\geq \phi_\sA^+(M)+\ell(1-d)}$ is independent of $n>\ell$.
Therefore, taking
$$
\widetilde{F}=\left(T_E^{\ell+1}(F)\right)_\sA^{\geq \phi_\sA^+(M)+\ell(1-d)} \quad \text{ and } \quad C_n=\left(T_E^{n}(F)\right)_\sA^{<\phi_\sA^+(M)+\ell(1-d)}
$$
would satisfy the desired property.
\end{proof}

\begin{cor}
\label{cor:SphericalEntropy}
Let $E$ be a $d$-spherical object in $\cD$, where $d\geq2$.
Suppose $E$ satisfies Conditions~\ref{cond:hasHomN0} and \ref{cond:existHeartPerp}. Then, the categorical entropy function of the spherical twist $T_E$ is given by:
$$
h_t(T_E)=
\begin{cases}
    (1-d)t & \text{ for } t\leq0, \\
    0 & \text{ for } t\geq0.
\end{cases}
$$
\end{cor}
\begin{proof}
This follows from Theorem~\ref{thm:Ouchispherical}, Lemma~\ref{lemma:SphericalMaxDegree}, and Proposition~\ref{prop:SphericalMaxDegree}.
\end{proof}

\begin{cor}[\cite{LiLiu}*{Theorem~A.2}]
\label{cor:LiLiuOuchi}
Consider a complex smooth projective variety $X$ of dimension $d\geq2$, and let $E\in\Db(X)$ be a $d$-spherical object with the associated spherical twist $T_E\in\Aut\Db(X)$.
Let $k$ be a nonzero integer, and $Y$ be a smooth projective variety with an exact equivalence $\Psi\colon\Db(X)\rightarrow\Db(Y)$. 

There does not exist a standard autoequivalence $\Phi\in\Aut_\std\Db(Y)$ such that $T_E^k=\Psi^{-1}\circ \Phi\circ\Psi$. In particular, $T_E^k$ is not conjugate to a standard autoequivalence of $X$.
\end{cor}
\begin{proof}
Recall that Conditions~\ref{cond:hasHomN0} and \ref{cond:existHeartPerp} are always satisfied for $\cD=\Db(X)$, see Remark~\ref{rmk:(c)(d)D(X)satisfy}.
Therefore, by Corollary~\ref{cor:SphericalEntropy}, $\tau^+(T_E)=0$ and $\tau^-(T_E)=1-d<0$. Then, according to \cite{FanFilip}*{Theorem~1.1(iii)}, we have
$$
\tau^+(T_E^k)=
\begin{cases}
    0  & \text{ if } k>0, \\
    -k(d-1) & \text{ if } k<0,
\end{cases}
\qquad
\tau^-(T_E^k)=
\begin{cases}
    -k(d-1)  & \text{ if } k>0, \\
    0 & \text{ if } k<0.
\end{cases}
$$
Thus, $\tau^+(T_E^k)>\tau^-(T_E^k)$ for any $k\neq0$. The desired statement then follows from Proposition~\ref{prop:EntropyBasic}\ref{propEntropy:conjugacyInvar} and Example~\ref{eg:stdautoeq}.
\end{proof}

\subsubsection{Entropy function of $\bP$-twists}
\label{subsubsec:P-twist}

In the following, we prove the analogues of Proposition~\ref{prop:SphericalMaxDegree}, Corollary~\ref{cor:SphericalEntropy}, and Corollary~\ref{cor:LiLiuOuchi} for $\bP$-twists. The arguments are essentially identical to their spherical twist counterparts, so we will only highlight the differences for the $\bP$-twists.

\begin{prop}
\label{prop:P-MaxDegree}
Let $E$ be a $\bP^d$-object satisfying Conditions~\ref{cond:hasHomN0} and \ref{cond:existHeartPerp}. Then, there exists a nonzero object $F\in\cD$ and a heart $\sA\subseteq\cD$ such that
$$
\lim_{n\rightarrow\infty}\frac{\phi_{\sA}^+(P_E^n(F))}{n}=0.
$$
In fact, there exists a nonzero object $\widetilde{F}$ and an integer $\ell>0$ such that for any $n>\ell$, there is an exact triangle
$$
\widetilde{F}\rightarrow P_E^n(F) \rightarrow C_n \xrightarrow{+1}
$$
with $\phi_{\sA}^-(\widetilde{F})>\phi_{\sA}^+(C_n)$.
Therefore, $\phi_\sA^+(P_E^n(F))$ stabilizes when $n>\ell$.
\end{prop}

\begin{proof}
Let $F$ be a nonzero object in $\cD$. Define
$$
M\coloneqq\Cone\left(\Hom^{\ast-2}(E,F)\otimes E\xrightarrow{h^\vee\cdot\id-\id\cdot h}\Hom^\ast(E,F)\otimes E\right)[1].
$$
Applying $P_E^{n-1}$ to the exact triangle
$$
F\rightarrow P_E(F)\rightarrow M \xrightarrow{+1},
$$
one obtains
$$
 P_E^{n-1}(F)\rightarrow P_E^n(F)\rightarrow M[-2d(n-1)]\xrightarrow{+1}.
$$
Then we have
$$
\xymatrix@C=.5em{
F \ar[rrrr]
&&&&
\ast \ar[rrrr] \ar[dll]
&&&&
\ast \ar[rr] \ar[dll]
&&
\cdots \ar[rr] 
&&
\ast \ar[rrrr]
&&&&
P_E^n(F) \ar[dll] 
\\
&&
M \ar@{-->}[ull]
&&&&
M[-2d] \ar@{-->}[ull]
&&&&&&&&
M[-2d(n-1)] \ar@{-->}[ull]
}
$$
Since $E$ satisfies Condition~\ref{cond:existHeartPerp}, there exists a heart $\sA\subseteq\cD$ such that (at least) one of the following statements holds:
\begin{itemize}
    \item there exists a nonzero object $A\in\sA$ such that $\Hom(E_{\sA}^{+},A)=0$, or
    \item $\phi_\sA^+(E)-\phi_\sA^-(E)\leq1$.
\end{itemize}

First, suppose there exists $A\in\sA$ such that $\Hom(E_{\sA}^{+},A)=0$.
Then, taking $F=A$ would satisfy the desired properties.
The proof is identical to the case of spherical twists, as we also have  $M_\sA^+\cong(E_\sA^+)^{\oplus p}$ in this scenario.

Second, suppose $\phi_\sA^+(E)-\phi_\sA^-(E)\leq1$.
By Lemma~\ref{lemma:concentrateInd-1}, there exists $F\in\cD$ such that $\Hom(E,F)\neq0$, and $\Hom(E,F[k])=0$ if $k<0$ or $k>2d-2$. Therefore, $N\coloneqq\RHom(E,F)\otimes E$ is nonzero and satisfies
$$
\phi_\sA^+(N)-\phi_\sA^-(N)\leq 2d-1.
$$
By the definitions of $M$ and $N$, there is an exact triangle
$$
N[1] \rightarrow M \rightarrow N \xrightarrow{+1}.
$$
Thus, we have
$$
\phi_\sA^+(M)-\phi_\sA^-(M)\leq 2d.
$$
Choose $\ell>0$ large enough so that
$$
\phi_\sA^-(F)>\phi_\sA^+(M[-2d\ell])=\phi_\sA^+(M)-2d\ell.
$$
Let $\widetilde{\phi}\coloneqq\phi_\sA^+(M)-2d\ell$.
Then we have the following picture, analogous to what we had in the case of spherical twists. The remainder of the proof follows identically.
\begin{center}
\begin{tikzpicture}
\draw (-1.5, 0) rectangle  ++(3,1) node[pos=.5, align=center] {$\cdots$};
\draw (1.5, 0) rectangle  ++(1.5,1) node[pos=.5, align=center] {$\sA[\widetilde{\phi}]$};
\draw (3, 0) rectangle  ++(2,1) node[pos=.5, align=center] {$\cdots$};
\draw (5, 0) rectangle  ++(3,1) node[pos=.5, align=center] {$\sA[\widetilde{\phi}-2d]$};
\draw (8, 0) rectangle  ++(2,1) node[pos=.5, align=center] {$\cdots$};
\draw (10, 0) rectangle  ++(3,1) node[pos=.5, align=center] {$\sA[\widetilde{\phi}-4d]$};
\draw (13, 0) rectangle  ++(2,1) node[pos=.5, align=center] {$\cdots$};
\filldraw[draw=black,fill=yellow] (-1.5,-1)rectangle ++(2.5,1) node[pos=.5, align=center] {$F$};
\filldraw[draw=black,fill=yellow] (-1.5,-2)rectangle ++(4.5,1) node[pos=.5, align=center] {$M[-2d(\ell-1)]$};
\filldraw[draw=black,fill=yellow] (1.5,-3)rectangle ++(6.5,1) node[pos=.5, align=center] {$M[-2d\ell]$};
\draw (2.25,-2.5) node {$\times$};
\filldraw[draw=black,fill=yellow] (5,-4)rectangle ++(8,1) node[pos=.5, align=center] {$M[-2d(\ell+1)]$};
\draw (6.5,-3.5) node {$\times$};
\filldraw[draw=black,fill=yellow] (10,-5)rectangle ++(5,1);
\draw (11.5,-4.5) node {$\times$};
\end{tikzpicture}
\end{center}
\end{proof}

\begin{cor}
\label{cor:PTwistEntropy}
Let $E$ be a $\bP^d$-object in $\cD$.
Suppose $E$ satisfies Conditions~\ref{cond:hasHomN0} and \ref{cond:existHeartPerp}. Then, the categorical entropy function of the $\bP$-twist $P_E$ is given by:
$$
h_t(P_E)=
\begin{cases}
    -2dt & \text{ for } t\leq0, \\
    0 & \text{ for } t\geq0.
\end{cases}
$$
\end{cor}
\begin{proof}
This follows from Theorem~\ref{thm:OuchiP-twistVersion}, Lemma~\ref{lemma:P-MaxDegree}, and Proposition~\ref{prop:P-MaxDegree}.
\end{proof}

\begin{cor}[= Theorem~\ref{thm:MainP-twist}]
Consider a complex smooth projective variety $X$ of dimension $2d$, and let $E\in\Db(X)$ be a $\bP^d$-object with the associated $\bP^d$-twist $P_E\in\Aut\Db(X)$.
Let $k$ be a nonzero integer, and $Y$ be a smooth projective variety with an exact equivalence $\Psi\colon\Db(X)\rightarrow\Db(Y)$. 

There does not exist a standard autoequivalence $\Phi\in\Aut_\std\Db(Y)$ such that $P_E^k=\Psi^{-1}\circ \Phi\circ\Psi$. In particular, $P_E^k$ is not conjugate to a standard autoequivalence of $X$.
\end{cor}
\begin{proof}
Recall that Conditions~\ref{cond:hasHomN0} and \ref{cond:existHeartPerp} are always satisfied for $\cD=\Db(X)$, see Remark~\ref{rmk:(c)(d)D(X)satisfy}.
Therefore, by Corollary~\ref{cor:PTwistEntropy}, $\tau^+(P_E)=0$ and $\tau^-(P_E)=-2d<0$. Then, according to \cite{FanFilip}*{Theorem~1.1(iii)}, we have
$$
\tau^+(P_E^k)=
\begin{cases}
    0  & \text{ if } k>0, \\
    -2dk & \text{ if } k<0,
\end{cases}
\qquad
\tau^-(P_E^k)=
\begin{cases}
    -2dk  & \text{ if } k>0, \\
    0 & \text{ if } k<0.
\end{cases}
$$
Thus, $\tau^+(P_E^k)>\tau^-(P_E^k)$ for any $k\neq0$. The desired statement then follows from Proposition~\ref{prop:EntropyBasic}\ref{propEntropy:conjugacyInvar} and Example~\ref{eg:stdautoeq}.
\end{proof}

\subsection{Polynomial entropy and Proof of Theorem~\ref{thm:PolyCat}}
\label{subsec:ProofOfThm2}

\begin{prop}
\label{prop:polycat}
Let $E$ be a $d$-spherical object in $\cD$, where $d\geq2$.
\begin{enumerate}[label=(\roman*)]
    \item\label{item:TwistPolyEnt-1} For $t<0$, $h_t^\pol(T_E)=0$.
    \item\label{item:TwistPolyEnt-2} If $E$ satisfies Conditions~\ref{cond:hasHomN0} and \ref{cond:existHeartPerp}, then $h_t^\pol(T_E)=0$ for $t>0$.
    \item\label{item:TwistPolyEnt-3} If $E$ satisfies Conditions~\ref{cond:hasHomN0} and \ref{cond:StabilityHeart}, then $h_0^\pol(T_E)=1$.
\end{enumerate}
The same statements hold for the $\bP$-twists.
\end{prop}

We will prove the proposition only in the case of spherical twists, as the proof for the $\bP$-twists is essentially identical.

\begin{proof}[Proof of \ref{item:TwistPolyEnt-1}]\renewcommand{\qedsymbol}{}
This is proved in \cite{FanFuOuchi}, so we omit the proof here.
\end{proof}

\begin{proof}[Proof of \ref{item:TwistPolyEnt-2}]\renewcommand{\qedsymbol}{}
Let $M=\RHom(E,G)\otimes E[1]$.
Applying $T_E^{n-1}$ to the exact triangle
$$
G\rightarrow T_E(G)\rightarrow M \xrightarrow{+1},
$$
one obtains
$$
 T_E^{n-1}(G)\rightarrow T_E^n(G)\rightarrow M[(1-d)(n-1)]\xrightarrow{+1}.
$$
Thus
\begin{align*}
\epsilon_t(G,T_E^n(G)) & \leq \epsilon_t(G,M) e^{(1-d)(n-1)t} + \epsilon_t(G,T_E^{n-1}(G)) \\
& \leq \epsilon_t(G,M) e^{(1-d)(n-1)t} + \epsilon_t(G,M) e^{(1-d)(n-2)t} + \epsilon_t(G,T_E^{n-2}(G)) \\
& \leq \cdots \\
& \leq \epsilon_t(G,M)\left(\sum_{k=0}^{n-1} e^{(1-d)kt}\right) +\epsilon_t(G,G).
\end{align*}
Let $t>0$. By Corollary~\ref{cor:SphericalEntropy}, $h_t(T_E)=0$ when $E$ satisfies Conditions~\ref{cond:hasHomN0} and \ref{cond:existHeartPerp}. Therefore,
\begin{align*}
h_t^\pol(T_E) & =\limsup_{n\rightarrow\infty}\frac{\log\epsilon_t(G,T_E^n(G))}{\log(n)} \\
& \leq \limsup_{n\rightarrow\infty}\frac{1}{\log(n)}\log \left(\epsilon_t(G,M)\left(\sum_{k=0}^{n-1} e^{(1-d)kt}\right) +\epsilon_t(G,G)\right) =0.
\end{align*}
To obtain the lower bound, recall from Proposition~\ref{prop:SphericalMaxDegree} that there exist nonzero objects $F,\widetilde{F}$ and an integer $\ell>0$ such that for any $n>\ell$, there is an exact triangle
$$
\widetilde{F}\rightarrow T_E^n(F) \rightarrow C_n \xrightarrow{+1}
$$
with $\phi_{\sA}^-(\widetilde{F})>\phi_{\sA}^+(C_n)$.
Observe that $\Hom(\widetilde{F},C_n)=\Hom(\widetilde{F},C_n[-1])=0$ for all $n>\ell$. Thus
$$
\Hom(\widetilde{F},T_E^n(F))\cong\Hom(\widetilde{F},\widetilde{F})\neq0.
$$
Let $G$ be a split generator. Then both $G\oplus F$ and $G\oplus\widetilde{F}$ are also split generators. We have
\begin{align*}
h_t^\pol(T_E) & =\limsup_{n\rightarrow\infty}\frac{\log\epsilon_t(G\oplus\widetilde{F},T_E^n(G\oplus F))}{\log(n)} \\
& \geq \limsup_{n\rightarrow\infty}\frac{\log\epsilon_t(\widetilde{F},T_E^n(F))}{\log(n)} \\
& \geq \limsup_{n\rightarrow\infty}\frac{\log(1)}{\log(n)}=0.
\end{align*}
This concludes the proof of \ref{item:TwistPolyEnt-2}.
\end{proof}

\begin{proof}[Proof of \ref{item:TwistPolyEnt-3}]
First, there is an upper bound
\begin{align*}
h_0^\pol(T_E) & =\limsup_{n\rightarrow\infty}\frac{\log\epsilon_0(G,T_E^n(G))}{\log(n)} \\
& \leq \limsup_{n\rightarrow\infty}\frac{1}{\log(n)}\log \left(\epsilon_0(G,M)\cdot n+\epsilon_0(G,G)\right)\leq 1.
\end{align*}
Therefore, it suffices to show that $h_0^\pol(T_E)\geq1$.

Let $\cA=\sP_\sigma(0,1]$ be the heart associated to the given Bridgeland stability condition on $\cD$ in Condition~\ref{cond:StabilityHeart}.
Since $E\in\sP_\sigma(-1,1]$, we have $\phi_\sA^+(E)-\phi_\sA^-(E)\leq1$.
By Lemma~\ref{lemma:concentrateInd-1}, there exists $F\in\cD$ such that $\Hom(E,F)\neq0$, and $\Hom(E,F[k])=0$ if $k<0$ or $k>d-2$. Therefore, $M=\RHom(E,F)\otimes E[1]$ is nonzero and satisfies
$$
\phi_\sA^+(M)-\phi_\sA^-(M)\leq d-1.
$$
Choose $\ell>0$ large enough so that
$$
\phi_\sA^-(F)>\phi_\sA^+(M[\ell(1-d)])=\phi_\sA^+(M)+\ell(1-d)\eqqcolon\widetilde{\phi}.
$$
Again, we have the following picture for $T_E^n(F)$.
\begin{center}
\begin{tikzpicture}
\draw (-1.5, 0) rectangle  ++(3,1) node[pos=.5, align=center] {$\cdots$};
\draw (1.5, 0) rectangle  ++(1.5,1) node[pos=.5, align=center] {$\sA[\widetilde{\phi}]$};
\draw (3, 0) rectangle  ++(2,1) node[pos=.5, align=center] {$\cdots$};
\draw (5, 0) rectangle  ++(3,1) node[pos=.5, align=center] {$\sA[\widetilde{\phi}+(1-d)]$};
\draw (8, 0) rectangle  ++(2,1) node[pos=.5, align=center] {$\cdots$};
\draw (10, 0) rectangle  ++(3,1) node[pos=.5, align=center] {$\sA[\widetilde{\phi}+2(1-d)]$};
\draw (13, 0) rectangle  ++(2,1) node[pos=.5, align=center] {$\cdots$};
\filldraw[draw=black,fill=yellow] (-1.5,-1)rectangle ++(2.5,1) node[pos=.5, align=center] {$F$};
\filldraw[draw=black,fill=yellow] (-1.5,-2)rectangle ++(4.5,1) node[pos=.5, align=center] {$M[(\ell-1)(1-d)]$};
\filldraw[draw=black,fill=yellow] (1.5,-3)rectangle ++(6.5,1) node[pos=.5, align=center] {$M[\ell(1-d)]$};
\draw (2.25,-2.5) node {$\times$};
\filldraw[draw=black,fill=yellow] (5,-4)rectangle ++(8,1) node[pos=.5, align=center] {$M[(\ell+1)(1-d)]$};
\draw (6.5,-3.5) node {$\times$};
\filldraw[draw=black,fill=yellow] (10,-5)rectangle ++(5,1);
\draw (11.5,-4.5) node {$\times$};
\end{tikzpicture}
\end{center}
Recall that $h_0(T_E)=0$ by Theorem~\ref{thm:Ouchispherical}\ref{item:Genki-1}. We claim that $h_{\sigma,0}(T_E)=0$. Indeed, by Lemma~\ref{lemma:entropyVSmassgrowth}\ref{item:entropyVSmassgrowth-1}, $h_{\sigma,0}(T_E)\leq h_0(T_E)=0$. On the other hand, by Remark~\ref{rmk:MassLowerBound}, there is a constant $C>0$ (arising from the support property of stability conditions) such that
$$
h_{\sigma,0}(T_E)=\limsup_{n\rightarrow\infty}\frac{1}{n}\log m_{\sigma,0}(\Phi^nG)\geq\limsup_{n\rightarrow\infty}\frac{\log C}{n}=0.
$$
This proves that $h_0(T_E)=h_{\sigma,0}(T_E)=0$.
By Lemma~\ref{lemma:entropyVSmassgrowth}\ref{item:entropyVSmassgrowth-2}, this implies that $h_0^\pol(T_E)\geq h_{\sigma,0}^\pol(T_E)$. Therefore, to obtain the desired lower bound for $h_0^\pol(T_E)$, it suffices to show that $h_{\sigma,0}^\pol(T_E)\geq1$.

From the above picture and Lemma~\ref{lemma:ExtremalDegree}\ref{item:NonzeroCohom}, one observes that $T_E^n(F)$ admits at least $n-\ell$ nonzero cohomology objects (with respect to $\sA$) when $n>\ell$, at degrees 
$$
\widetilde{\phi},\ \widetilde{\phi}+(1-d),\ \ldots,\ \widetilde{\phi}+(n-\ell-1)(1-d).
$$
Therefore, by Remark~\ref{rmk:MassLowerBound}, we have $m_{\sigma,0}(T_E^n(F))>C(n-\ell)$. Thus,
$$
h_{\sigma,0}^\pol(T_E)\geq\limsup_{n\rightarrow\infty}\frac{\log m_{\sigma,0}(T_E^n(F))}{\log(n)}
\geq\limsup_{n\rightarrow\infty}\frac{\log(C(n-\ell))}{\log(n)}=1.
$$
This completes the proof.
\end{proof}

\begin{cor}[= Theorem~\ref{thm:PolyCat}]
Let $X$ be a complex smooth projective variety of dimension $d\geq2$. Let $E\in\Db(X)$ be a spherical object (resp.~$\bP$-object) with the associated spherical twist $T_E\in\Aut\Db(X)$ (resp.~$\bP$-twist $P_E\in\Aut\Db(X)$). Then
\begin{enumerate}[label=(\roman*)]
    \item For $t\neq0$, $h_t^\pol(T_E)=0$ (resp.~$h_t^\pol(P_E)=0$).
    \item If there exists a Bridgeland stability condition $\sigma=(Z,\sP)$ on $\cD$ such that $E\in\sP(-1,1]$, then $h_0^\pol(T_E)=1$ (resp.~$h_0^\pol(P_E)=1$).
\end{enumerate}
\end{cor}
\begin{proof}
This follows from Proposition~\ref{prop:polycat}, and the fact that Conditions~\ref{cond:hasHomN0} and \ref{cond:existHeartPerp} are always satisfied for $\cD=\Db(X)$, cf.~Remark~\ref{rmk:(c)(d)D(X)satisfy}.
\end{proof}

\bigskip
\bibliography{ref}
\bibliographystyle{alpha}

\ \\

\noindent Yu-Wei Fan \\
\textsc{Yau Mathematical Sciences Center, Tsinghua University\\
Beijing 100084, China}\\
\texttt{ywfan@mail.tsinghua.edu.cn}

\end{document}